\documentclass[a4paper, 11pt, parskip]{scrartcl}
 
\usepackage[T1]{fontenc}
\usepackage{amssymb, amsfonts, amsmath, amsthm, dsfont}
\usepackage[english]{babel}
\usepackage{ae, aecompl} 
\usepackage{url} 
\usepackage{color}
\usepackage{arydshln} 
\usepackage{fancyhdr}
\usepackage{rotating}
\usepackage{nicefrac}
\usepackage{subfigure}
\usepackage{makeidx}
\usepackage[all]{xy}

\usepackage[pdfauthor={Oliver Braun and Renaud Coulangeon},
            pdftitle={Perfect Lattices for Imaginary Quadratic Number Fields}]{hyperref}
\theoremstyle{plain}
\newtheorem{thm}{Theorem}[section]
\newtheorem{lemma}[thm]{Lemma}

\theoremstyle{definition}
\newtheorem{definition}[thm]{Definition}

\newtheorem{remark}[thm]{Remark}

\newtheorem{example}[thm]{Example}
\newtheorem{algorithm}[thm]{Algorithm}

 \newcommand{\MZ}{\mathds{Z}} \newcommand{\MQ}{\mathds{Q}}
\newcommand{\MR}{\mathds{R}} \newcommand{\MC}{\mathds{C}}

\newcommand{\beweisende}{\hspace*{\fill} $\square$ }
 
\newcommand{\GL}{\operatorname{GL}}

\newcommand{\PSL}{\operatorname{PSL}}

\renewcommand{\O}{\mathcal{O}} \newcommand{\V}{\mathcal{V}}  
\newcommand{\stab}{\operatorname{Stab}} 

\renewcommand{\a}{\mathfrak{a}}  \renewcommand{\c}{\mathfrak{c}}
\newcommand{\cl}{\mathcal{C}\ell} 
\newcommand{\trace}{\mathrm{Trace}} \renewcommand{\H}{\mathcal{H}} \newcommand{\A}{\mathcal{A}} \newcommand{\B}{\mathcal{B}}

\newcommand{\st}{\mathrm{St}}
\newcommand{\tr}{\mathrm{Tr}} 
\newcommand{\p}{\mathfrak{p}}

\title{Perfect Lattices for Imaginary Quadratic Number Fields}
\date{}
\author{Oliver Braun \\
Lehrstuhl D f\"{u}r Mathematik, RWTH Aachen University,\\
Templergraben 64, D-52062 Aachen, Germany\\
\url{oliver.braun1@rwth-aachen.de}\\
\mbox{} \\
Renaud Coulangeon \\
Univ. Bordeaux, IMB, UMR 5251, F-33400 Talence, France.\\
CNRS, IMB, UMR 5251, F-33400 Talence, France.\\
\url{Renaud.Coulangeon@math.u-bordeaux1.fr}}

\begin{document}

\maketitle
\begin{abstract}
We present an adaptation of Voronoi theory for imaginary quadratic number fields of class number greater than 1. This includes a characterisation of extreme Hermitian forms which is analogous to the classic characterisation of extreme quadratic forms as well as a version of Voronoi's famous algorithm which may be used to enumerate all perfect Hermitian forms for a given imaginary quadratic number field in dimensions 2 and 3. We also present an application of the algorithm which allows to determine generators of the general linear group of an $\O_K$-lattice.
\end{abstract}

\section{Introduction}
The notion of perfect lattices (resp. quadratic form) first appeared in Voronoi's celebrated paper \cite{Vo1,Vo2}, as part of the characterization of the lattices corresponding to locally densest regular sphere packings, the so-called \emph{extreme lattices}. As pointed out in \cite[p.105]{mar03}, the concept is already  visible in Korkine and Zolotareff's paper \cite{KZ77}, as well as its relevance to packing density, although with no name associated to it. Maybe as importantly, Voronoi proved that there are finitely many perfect forms in a given dimension, up to integral equivalence, and described an algorithm to compute them. This part of Voronoi's theory has had numerous consequences. Indeed, Voronoi's algorithm not only computes perfect forms but builds a cellular complex that is acted on by $\GL_n(\MZ)$. This observation is at the core of the (co)-homology computations for arithmetic groups developed by Ash \cite{Ash84}, Soul\'e \cite{Sou78} and others (a very comprehensive account on this beautiful topic can be found in Appendix A of \cite{Ste07} written by Paul E. Gunnells).

An extension of these concepts to lattices over number fields is underlying Humbert's article \cite{Hum49}, although no notion of perfection is considered there. Later, Koecher \cite{Koe60,Koe61} undertook a vast generalization of Voronoi's theory which encompasses the case of positive definite quadratic forms (resp. lattices) over number fields. More recently, Watanabe \cite{Wat00,Wat03} built a very general framework for Hermite-like functions over algebraic groups, including the classical Hermite function of Euclidean lattices, viewed as a function on $\GL_n(\mathds A_{\MQ})$, the adelic general linear group. When applied to $\GL_n(\mathds A_{K})$, $K$ a number field, Watanabe's theory gives rise to a theory of \emph{lattices over number fields} which is in general substantially different from that of Koecher, but coincides with it in the case of imaginary quadratic fields, to which we restrict in this paper. Watanabe's theory is formulated in an adelic language, but for the general linear group over a number field, it can be translated in terms of so-called \emph{Humbert forms} and projective modules over Dedekind domains (see \cite{cou04}). This is the approach we adopt here (see section \ref{preliminaries} for details). While the adelic setting provides much more uniform formulations and statements (for instance, there is no need for a distinction between free and non-free lattices with the adelic point of view), it is certainly not well-suited for computations, which are the main concern of this paper. 

To be more precise, we study the notion of Hermitian perfect forms over an imaginary quadratic field $K=\MQ(\sqrt{-d})$, where $d$ is a square free positive integer, and give complete classification results in dimension $2$ and $3$ over fields of small discriminant. The main technical ingredient is an adaptation of the Voronoi algorithm which we explain in section \ref{algorithm}. One main difference with the classical situation (over $\MZ$), is that one has to take into account the non-principality of the ring of integers $\O_K$. This leads to define one Hermite function for each \emph{type} of $\mathcal O_K$-lattice, which are in one-to-one correspondence with the ideal classes of $K$ (see section \ref{preliminaries}), and run the algorithm separatly for each of those. The computational results are presented in section \ref{results}. A more detailed account, including the code used to perform the computations with MAGMA \cite{magma} is available at \url{http://www.math.rwth-aachen.de/~Oliver.Braun/perfect.html}
We have to mention that similar computations have already been carried out by Yasaki in \cite{Yas10}, but only for free lattices, and with a slightly different definition of the Hermite function.

In section \ref{opgen}, following a general strategy described by Opgenorth in \cite{Opg01}, we show how to use Voronoi graphs of perfect hermitian forms to compute generators of groups of the type $G=\GL (\O_K^{n-1} \oplus \a)$. When $n=2$, this includes classical Bianchi groups $\GL_2(\mathcal O_K)$, for which a general procedure 
is known, a description of which can be found in the book by Elstrodt, Gr\"unewald and Mennicke (see \cite [chap. 7]{EGM98}). Note that the procedure described there presupposes an explicit description of a fundamental domain for the action of $G$ on the three-dimensional upper half-space $\mathds H = \mathds C \times \left] 0, \infty\right[ $, which can be obtained in general using an algorithm due to Swan  \cite{Swa71} (see \cite{RF11} for examples of implementation). In comparison, our method does not rely on any \textit{a priori} description of a fundamental domain, and can be applied to any dimension $n$. It would thus theoretically have a wider scope.

Most of the results presented in this paper originate from the first author's bachelor thesis, written under the supervision of Professor Gabriele Nebe at RWTH Aachen University. The first author would like to thank her for her support and would also like to acknowledge the support of Bertrand Meyer during the completion of the thesis.

\section{Preliminaries}\label{preliminaries}

Throughout this paper, $K/\MQ$ will be an imaginary quadratic number field. For the remainder of this article we shall fix an embedding $K\hookrightarrow \MC$ so that, by abuse of notation, we can write $K\subseteq \MC$.
\\
We denote by $\O_K$ the integral closure of $\MZ$ in $K$, by $\cl_K$ the ideal class group of $K$, of cardinality $h_K$, and we let $\a_1, ... , \a_{h_K}$ be a set of representatives of the ideal classes which are chosen to be integral and of minimal norm.

By $\MZ^n, ~ \MQ^n , ~ K^n$ etc. we mean sets of row vectors with $n$ entries. As usual, we denote by $A^*$ the conjugate transpose of a matrix with complex entries, a square matrix being Hermitian if it equals its conjugate transpose. The set $\H_n$ of Hermitian $n\times n$-matrices, equipped with the scalar product
\begin{equation*}\H_n \times \H_n \to \MR, ~ (\A,\B)\mapsto \trace ( \A\B ) \end{equation*}
is a Euclidean vector space of dimension $n^2$ over $\MR$, which we identify with the set of Hermitian forms on $\MC^n$. The cone of positive definite Hermitian forms then corresponds to the cone $\H_n^+$ of positive definite Hermitian matrices in $\H_n$.

For $\A\in\H_n$ and $x\in K^n$ we set $\A[x]:=x\A x^*$.

\begin{definition}A lattice $L$ in $K^n$ is a finitely generated $\O_K$-submodule of $K^n$ such that $K\otimes_{\O_K} L \cong K^n$. We set $\GL(L):=\{\varphi \in \GL_n(K) ~|~ \varphi(L) \subseteq L\}$. Two $\O_K$-lattices $L$ and $M$ are \emph{isomorphic} if there exists $\varphi \in \GL_n(K)$ such that $\varphi(L)=M$.\end{definition}
 
\begin{thm}[Steinitz] Any $\O_K$-lattice $L$ in $K^n$ is isomorphic to a direct sum of fractional ideals of $K$, $L\cong \mathfrak{c}_1 \oplus ... \oplus \mathfrak{c}_n$, and the ideal class $[\mathfrak{c}_1 \cdot ... \cdot \mathfrak{c}_n] \in \cl_K$ is a complete invariant of the isomorphism class of $L$. The ideal class $[\mathfrak{c}_1 \cdot ... \cdot \mathfrak{c}_n] =: \st(L)$ is called the Steinitz Class of $L$.\end{thm}
\textbf{Proof:} \cite[Theorems 81:3, 81:11]{ome}\beweisende

In other words, any $\O_K$-lattice $L$ in $K^n$ may be written as $\mathfrak{c}_1 e_1\oplus ... \oplus \mathfrak{c}_ne_n$ where the $\mathfrak{c}_i$s are fractional ideals and $\left\lbrace e_1, \dots e_n\right\rbrace $ is a $K$-basis of $K^n$, and two lattices $\bigoplus \mathfrak{c}_i e_i$ and $\bigoplus \mathfrak{d}_i f_i$ are isomorphic if and only if $[\mathfrak{c}_1 \cdot ... \cdot \mathfrak{c}_n]=[\mathfrak{d}_1 \cdot ... \cdot \mathfrak{d}_n]$.

Note that a lattice $\mathfrak{c}_1 \oplus ... \oplus \mathfrak{c}_n$ is isomorphic to a lattice of the form $\O_K^{n-1}\oplus \mathfrak{c}_1 \cdot ... \cdot \mathfrak{c}_n$, so $\O_K^{n-1}\oplus \a_i$ is a complete set of representatives of the isomorphism classes of $\O_K$-lattices, if $i$ ranges from $1$ to $h_K$.

For a lattice $L=\mathfrak{c}_1 e_1\oplus ... \oplus \mathfrak{c}_ne_n$ in $K^n$ and a Hermitian form $\A$ on $K^n$ we have the following definitions which are cited from \cite{cou04} \footnote{Note that, for simplicity, we have omitted the squares which are present in the original definition in \cite{cou04}.}: for $x=\sum_{i=1}^n x_i e_i \in L$ we define the integral ideal $\a_x := x_1 \mathfrak{c}_1^{-1} + ... + x_n \mathfrak{c}_n^{-1}$, which may be interpreted as a "g.c.d. of the coefficients of $x$". 
Then we define the minimum, determinant and Hermite invariant of $\A$ with respect to $L$ respectively as
\begin{equation}\label{min} \min\nolimits _L (\A) := \min_{x\in L-\{0\}} \frac{\A[x]}{N(\a_x)}, \end{equation}

\begin{equation*}\label{det} \det\nolimits _L(\A):=N(\c_1 \cdot ... \cdot \c_n) \det \A \end{equation*}
and
\begin{equation*} \gamma_L(\A):= \frac{ \min\nolimits _{L} (\A)  } { ( \det \nolimits _L (\A) )^{1/n} }. \end{equation*}
Finally, we set
\begin{equation*}\gamma_L:= \sup_{\A \in \H_n^+} \gamma_L (\A) \end{equation*}

\begin{remark}Note that for $\A\in\H_n^+$, $\alpha \in \MR_{>0}$ and $U\in\GL(L)$ we have
\begin{equation*}\gamma_L(\A)=\gamma_L(\alpha \A), \qquad \gamma_L(\A) = \gamma_L(\A[U]) \end{equation*}
so that $\gamma_L$ defines a function on $\MR_{>0}\backslash \H_n^+ / \GL(L)$.\end{remark}

\begin{remark} An important feature of formula (\ref{min}), to be used repeatedly, is the fact that for any $x \in L-\{0\}$ and $\alpha \in K^*$, one has
\begin{equation*}
 \frac{\A[x]}{N(\a_{x})} = \frac{\A[\alpha x]}{N(\alpha \a_{x})} =\frac{\A[\alpha x]}{N(\a_{\alpha x})}. 
\end{equation*}
Consequently, for a given form $\A$, the map $x \mapsto \frac{\A[x]}{N(\a_{x})}$ may be viewed as a function on the projective space $\mathbb P(K^n)$. Also, when $n=2$, the $3$-dimensional hyperbolic space $\mathbb H:=\MC \times \left] 0,\infty\right[$, acted on by $\PSL_2(\O_K) \subset \PSL_2(\MC)$ provides a geometric model for the cone $\H_2^+/\PSL_2(\O_K)$ that has been studied by many authors (see \cite{men,flo83,RF11} and \cite[chapter 7]{EGM98}). In this setting, the quantity $\frac{\A[x]}{N(\a_{x})}$ measures the hyperbolic distance from the point $\A \in \H_2^+$ to the cusp associated to $x$. 
\end{remark}

The following lemma is easily verified.

\begin{lemma}Let $L, ~ L'$ be two $n$-dimensional $\O_K$-lattices with $\st(L)=\st(L')$. Then $\gamma_L = \gamma_{L'}$.\end{lemma}

This lemma shows that it is sufficient to calculate the constants $\gamma_i := \gamma_{L_i}$ for the lattices $L_i := \O_K^{n-1} \oplus \a_i$ in order to obtain
\begin{equation*} \gamma_{n,K}:=  \max _{1\leq i \leq h_K} \gamma_i \end{equation*}
which we call the ($n$-dimensional) Hermite Constant of the field $K$.

Note that in the case $K=\MQ$ this definition of Hermite's constant coincides with the traditional definition from the theory of $\MZ$-lattices.

In order to determine $\min\nolimits _ {L_i} ( \A)$ we can use the following lemma.

\begin{lemma} \begin{equation*} \min\nolimits_{L_i}(\A) = \min_{1\leq j \leq h_K} \min_{\substack{ x\in L_i - \{0\} \\ [\a_x] = [\a_j] }} \frac{ \A[x] } {N(\a_x)} = \min_{1\leq j \leq h_K} \min_ {\substack{ x \in L_i - \{0\} \\ \a_x = \a_j }} \frac{\A[x]}{N(\a_x)} \end{equation*} \end{lemma}
\textbf{Proof:} The first equality is clear. As for the second one, we note that for any $x\in L_i - \{0\}$,  there exist $j \in \left\lbrace 1, \dots, h_K\right\rbrace $ and $\alpha \in K^*$ such that $\a _j = \alpha \a_x = \a_{\alpha x}$, and we can conclude using the equality $\frac{\A[x]}{N(\a_x)}=\frac{\A[\alpha x]}{N(\a_{\alpha x})}$. 
\beweisende

\section{Perfection and Eutaxy}

In this section we will define appropriate generalisations of the concepts of perfection and eutaxy from the Voronoi theory of $\MZ$-lattices so that we can characterize extreme Hermitian forms.

In order to develop this characterization we need a suitable notion of shortest vectors of a Hermitian form $\A$. Since $\frac{\A[x]}{N(\a_x)} = \frac{ \A [\alpha x]}{ N(\a_{\alpha x}) }$ for all $x\in L$ and all $0\neq \alpha \in \O_K$ (even for all $\alpha \in K^*$), rather than vectors in $L$, we will consider their images in the projective space $\mathbb P(K^n)$. We will make use of the following lemma :
\begin{lemma}For any positive constant $C$, the set
\begin{equation*}  \left\{ x \in L_i - \{0\} ~|~ \frac{\A[x]}{N(\a_x)} = C\right\} / \sim \end{equation*}
is finite, where the equivalence relation $\sim$ on $L_i$ is defined by
\begin{equation*} v\sim w :\Longleftrightarrow ~\exists ~ \kappa \in K^* ~:~ v = \kappa w \end{equation*}
\end{lemma}
\textbf{Proof:} There is an injective map from the set $\{ x \in L_i - \{0\} ~|~ \frac{\A[x]}{N(\a_x)} = C \} / \sim$ into the finite set $\{x \in L_i ~|~ \A[x] \leq C\tilde{m} \}$ where $\tilde{m} := \max_{1\leq j \leq h_K} N(\a_j)$.\beweisende
\begin{definition}
The set of minimal vectors of $\A$ with respect to $L_i$ is  the finite set \begin{equation*}  S_i(\A) :=\left\{ x \in L_i - \{0\} ~|~ \frac{\A[x]}{N(\a_x)} = \min\nolimits _{L_i}(\A) \right\} / \sim \end{equation*}

\end{definition}
\begin{definition}We call $\A$ $\gamma_i$-perfect if
\begin{equation*}\langle x^* x ~|~ x \in S_i(\A) \rangle _\MR = \H_n, \end{equation*}
that is if the matrices $x^*x$ , as $x$ runs through $S_i(\A)$, span the whole $\MR$-vector space of Hermitian $n\times n$-matrices.\\
We call $\A$ $\gamma_i$-eutactic if for every $x\in S_i(\A)$ there is a $\lambda_x \in \MR_{>0}$ such that
\begin{equation*} \overline{\A^{-1}} = \sum_{x\in S_i(\A)} \lambda_x x^*x . \end{equation*}
We say that $\A$ is $\gamma_i$-extreme if it is a local maximum of the function $\gamma_{L_i}$.
\end{definition}

Note that $\gamma_i$-perfection is a property of the equivalence classes in $\MR_{>0} \backslash \H_n^+ / \GL(L_i)$.

As is the case in classical Voronoi theory of $\MZ$-lattices, $\gamma_i$-perfection entails the following properties which may be proven analogously, see for example \cite{mar03}.
\begin{thm}If $\A\in\H_n^+$ is $\gamma_i$-perfect then it is uniquely determined by $\min\nolimits _{L_i}(\A)$ and $S_i(\A)$.\end{thm}
\begin{remark}For any $\gamma_i$-perfect $\A\in\H_n^+$ we have $|S_i(\A)|\geq \dim_\MR  \H_n =n^2$ and $\langle S_i(\A) \rangle_\MR = \MC^n$.\end{remark}

The finiteness, up to $\GL(L_i)$-equivalence, of the set of $\gamma_i$-perfect forms, which in the case $K=\MQ$ was first proved by Voronoi in \cite{Vo1}, will turn out to be crucial. In our context, the precise statement, for which we may refer to \cite{meyer09}, is as follows

\begin{thm}[\cite{meyer09}]\label{finite} The set
\begin{equation*}\{ [\A] \in \MR_{>0}\backslash \H_n^+ / \GL(L_i) ~|~ \A ~ \gamma_i \text{-perfect}\} \end{equation*}
is finite.
\end{thm}

\textbf{Proof:} this is a particular case of \cite[Proposition 3.10 (1)]{meyer09}, the proof of which essentially relies on \emph{Humbert reduction theory} \cite{Hum49}, an alternative to Hermite reduction theory in the more general setting of $\O_K$-lattices.
\beweisende

As in classical Voronoi theory over $\MZ$ we have the following result which may be proved by adapting the methods from \cite[Chapter 3]{mar03}. Again, this can also be viewed as a particular case of Meyer's general result \cite[Theorem 3.9]{meyer09}

\begin{thm}\label{voronoitheorem}$\A$ is $\gamma_i$-extreme if and only if it is $\gamma_i$-perfect and $\gamma_i$-eutactic.\end{thm}

The following theorem will provide a way to simplify the calculation of the Hermite constant $\gamma_{n,K}$.

\begin{thm}\label{simplification}Let $L=\bigoplus_{j=1}^n \c_j$ be an $\O_K$-lattice and $\p$ any fractional ideal of $K$. Furthermore let $\overline{ \phantom{X} }$ be the non-trivial Galois automorphism of $K/\MQ$. Then for all $\A\in\H_n^+$ we have
\begin{equation*}\gamma_{L}(\A) = \gamma_{\p L}(\A) \quad \text{ and } \quad \gamma_L (\A) = \gamma_{\overline{L}}(\overline{\A}). \end{equation*}
Then the following are equivalent :
\begin{enumerate}
\item{$\A$ is perfect over $L$}
\item{$\A$ is perfect over $\p L$}
\item{$\overline{\A}$ is perfect over $\overline{L}$.}
\end{enumerate}
\end{thm}

\textbf{Proof:} Consider an arbitrary $\A \in \H_n^+$. We will show $\gamma_{L}(\A) = \gamma_{\p L}(\A)$. Let $x=\sum_{j=1}^n x_j e_j \in \p L$, we then have
\begin{equation*} \frac{ \A[x] } { N(x_1 (\p \c_1)^{-1} +...+ x_n (\p \c_n)^{-1}) } = N(\p) \frac{\A[x]} {N( x_1 \c_1^{-1} + ... + x_n \c_n^{-1} ) } = N(\p)\frac{\A[x]}{N(\a_x)} \end{equation*}
which implies $\min\nolimits _{\p L } (\A) = N(\p) \min\nolimits _L (\A)$ since the value $\frac{ \A[x] }{ N(\a_x) }$ only depends on the "direction" of $x$, i.e. the $K$-subspace it generates, and because any lattice $L$ in $K^n$ contains vectors of any given direction.

Furthermore $\det\nolimits_{\p L}(\A) = N(\p)^n \det \nolimits _L (\A)$, from which we obtain $\gamma_{\p L}(\A) = \gamma_{L}(\A)$.

It is also easily verified that we have $\gamma_{\overline{L}}(\overline{\A}) = \gamma_L (\A)$.

Theorem \ref{voronoitheorem} now implies that perfection of $\A$ over $L$ and $\p L$ is equivalent and that this is also equivalent to perefection of $\overline \A$ over $\overline L$ because the local maxima of Hermite's function are perfect Hermitian forms.\beweisende

Note that all significant invariants of $\A$ over $L$ and $\p L$ and of $\overline \A$ over $\overline L$ are identical because of the equality $\GL(L) = \GL(\p L)$ and the isomorphism $\GL(L) \cong \GL(\overline L)$.\\
Furthermore we have the equality
\begin{equation*}\st(\p L) = [\p ]^n \st(L)\end{equation*}
so that in order to determine $\gamma_{n,K}$, we need only to consider the lattices $\O_K^{n-1} \oplus \mathfrak{c}$, where $\mathfrak{c}$ runs over a set of representatives of $\mathrm{Gal}(K/\MQ)\backslash \cl_K / \cl_K^n$.

The following example illustrates the simplification in the process of determining $\gamma_{n,K}$.

\begin{example}Let $K/\MQ$ be a number field with cyclic ideal class group, $\cl_K \cong C_{h_K}$. Then the subgroup of all $n^{\text{th}}$ powers in $\cl_K$ is of order $\frac{h_K}{\mathrm{gcd}(h_K,n)}$, so that the quotient group $\cl_K/\cl_K^n$ is cylic of order $\mathrm{gcd}(h_K,n)$. Therefore, in this case, in order to calculate $\gamma_{n,K}$ instead of $h_K$ lattices we need only to consider $\mathrm{gcd}(h_K,n)$ lattices.\end{example}

\section{The Voronoi-algorithm}\label{algorithm}

In this section we will present an adaptation of Voronoi's algorithm to the situation of an imaginary quadratic number field. As in the classical case the algorithm enumerates all $\gamma_i$-perfect Hermitian forms. In addition we obtain a face-to-face tesselation of the cone of positive definite Hermitian forms $\H_n^+$. 

\begin{definition}Let $\A\in\H_n$ be positive definite. The $\gamma_i$-Voronoi domain of $\A$ is the convex closed cone
\begin{equation*}\V_i(\A) := \left\{ \sum_{x\in S_i(\A)} \lambda_x \cdot x^*x ~|~ \lambda_x \in \MR_{\geq 0} \right\}\subset \overline{\H_n^+}. \end{equation*}
where $\overline{\H_n^+}$ is the cone of positive semi definite matrices.
\end{definition}

\begin{remark}We note that $\A \in \H_n^+$ is $\gamma_i$-perfect if and only if $\V_i(\A)$ has non-empty interior (in $\H_n$), or equivalently, if and only if $\V_i(\A)$ is not contained in any hyperplane of $\H_n$. Note also that the \emph{relative} interior of $\V_i(\A)$ is $\displaystyle \left\{ \sum_{x\in S_i(\A)} \lambda_x \cdot x^*x ~|~ \lambda_x \in \MR_{> 0} \right\}$.\end{remark}

In what follows, we use the term \emph{facet} to mean an $(n-1)$-dimensional face of a polyhedral cone.
\begin{definition}Let $\A\in\H_n^+$ be $\gamma_i$-perfect and let $\mathcal{S}$ be a facet of the Voronoi domain $\V_i(\A)$. Then we call $0\neq R \in \H_n$ a facet vector if
\begin{equation*}\trace(RS)=0 ~\forall ~ S\in\mathcal{S}, \qquad \trace(RT)\geq 0 ~\forall ~ T\in \V_i(\A). \end{equation*}
\end{definition}
In other words, $R$ is a facet vector of the facet $\mathcal{S}$ if and only if the following two conditions are satisfied :
\begin{enumerate}
\item{$\trace (x^*xR)=R[x]=0$ for all $x\in S_i(\A)$ satisfying $x^*x \in \mathcal{S}$}
\item{$\trace (x^*xR)=R[x]>0$ for all $x\in S_i(\A)$ satisfying $x^*x \notin \mathcal{S}$.}
\end{enumerate}
This also shows that a facet vector can be chosen as the solution of a homogeneous system of linear equations with coefficients in $K$.

\begin{thm}For $n\geq 2$ let $\A\in\H_n$ be $\gamma_i$-perfect and let $\mathcal{S}$ be a facet of $\V_i(\A)$ with facet vector $R$. Then $R$ is indefinite and there is a lattice vector $x\in L_i$ satisfying $R[x]<0$.\end{thm}
\textbf{Proof:} Assume that $R$ is not indefinite. Since $R[x]\geq 0$ for all $x\in \MC^n$ we have
\begin{equation}\label{equality}\{v\in\MC^n ~|~ R[v]=0\} = \{v\in\MC^n ~|~ vR =0\}. \end{equation}
The righthand side of \eqref{equality} is a $d$-dimensional subspace of $\MC^n$ for some $d<n$. This implies that $\langle v^*v ~|~ \trace(v^*vR)=0\rangle \leq \H_n$ is a subspace of dimension at most 
\begin{equation*}d^2 \leq (n-1)^2 < n^2-1 = \dim(\H_n)-1\end{equation*}
However, $\mathcal{S}$ is a subspace of dimension $n^2-1$ which is generated by $v^*v$, where $v$ runs through  a set of isotropic vectors of $R$. This is a contradiction, therefore $R$ is indefinite. Consequently, there exists $x\in\MC^n$ such that $R[x]<0$. Since the map $x\mapsto R[x]$ is continuous and $K^n \subseteq \MC^n$ is a dense subset, we obtain a vector $y\in K^n$ satisfying $R[y]<0$. So we can also find a lattice vector in $L_i$ with the same property.\beweisende

\begin{thm}\label{interior}For $\A_1, ~\A_2 \in \H_n^+$ the following statements hold.
\begin{enumerate}
\item{If $T\in \V_i(\A_2)$ is contained in the relative interior of $\V_i(\A_1)$ we have $\V_i(\A_1) \subseteq \V_i(\A_2)$.}
\item{If a Hermitian form is contained within the interior of the $\gamma_i$-Voronoi domain of a perfect form, then it cannot be contained in the $\gamma_i$-Voronoi domain of any other perfect form.}
\end{enumerate}
\end{thm}
\textbf{Proof:} This is a simple adaptation of \cite[Theorem 7.1.12]{mar03}. \beweisende

Voronoi's algorithm is mainly based on the following theorem.
\begin{thm}\label{contiguous}Let $\A\in\H_n^+$ be $\gamma_i$-perfect and let $R$ be a facet vector associated with the facet $\mathcal{S}$ of $\V_i(\A)$. Put $S:=\{x\in S_i(\A) ~|~ x^*x \in \mathcal{S}\}$, $m:=\min\nolimits _{L_i}(\A)$.\\
For $t\in \MR$ define
\begin{equation*}\A_t := \A+t\cdot R \in \H_n.\end{equation*}
Then there is exactly one $\rho\in\MR_{>0}$ such that the following properties hold.
\begin{enumerate}
\item{If $0<t<\rho$ then $\A_t$ is not $\gamma_i$-perfect and $\min\nolimits_{L_i}(\A_t)=m$.
If $t>\rho$ then either $\A_t$ is not positive definite or $\min\nolimits_{L_i}(\A_t)<m$.}
\item{$S_i(\A_t)=S$ if $0<t<\rho$.}
\item{For negative values of $t$ $\A_t$ is not positive definite or $\min\nolimits_{L_i}(\A_t)<m$.}
\item{$\A_\rho$ is $\gamma_i$-perfect and $\min\nolimits_{L_i}(\A_\rho)=m$. $\mathcal{S}=\V_i(\A)\cap \V_i(\A_\rho)$. $\A$ and $\A_\rho$ are the only $\gamma_i$-perfect forms whose Voronoi domains contain $\mathcal{S}$.}
\end{enumerate}
\end{thm}
\textbf{Proof:} There exists $x\in L_i$ such that $R[x]<0$. So $\A_t$ is indefinite provided that $t$ is sufficiently large. Put
\begin{equation*}\rho:=\inf\{ t>0 ~|~ \min\nolimits_{L_i}(\A_t)<m \text{ or } \A_t \text{ is not positive definite}\} \end{equation*}
Note that $\rho>0$.

Now let $0<t<\rho$. Then we have $\min\nolimits_{L_i}(\A_t) \geq m$ and since $R[x]=0$ for all $x\in S$ equality holds. Therefore we can also conclude $S\subseteq S_i(\A_t)$.\\
For $y\in S_i(\A_t)$ assume $R[y]<0$. Then for all $t'>t$
\begin{equation*}\frac{\A_{t'}[y]}{N(\a_y)} < \min\nolimits_{L_i}(\A_t) = \min\nolimits_{L_i}(\A)\end{equation*}
Analogously we obtain this inequality for all $t'<t$ if we assume that $R[y]>0$. In conclusion we have $R[y]=0$ and therefore $y\in S$.\\
Since $\{x^*x ~|~ x\in S\}$ generates a hyperplane in $\H_n$ we see that $\A_t$ is not $\gamma_i$-perfect, which proves assertions 1 and 2.

In order to prove 3. consider $y\in S_i(\A)-S$ for $t<0$. Then we have $R[y]>0$ and
\begin{equation*} \frac{\A_t[y]}{N(\a_y)} = \frac{\A[y]}{N(\a_y)} + \frac{tR[y]}{N(\a_y)} < m. \end{equation*}

With regard to the fourth statement choose $y\in S_i(\A_\rho)$ such that $R[y]<0$. Such a choice is possible because $\min\nolimits_{L_i}(\A_t) < \min\nolimits_{L_i} (\A)$ for all $t>\rho$. In this situation $\langle x^*x ~|~ x\in S \cup \{y\} \rangle _\MR  = \H_n$ and $\A_\rho$ is $\gamma_i$-perfect.\\
Now consider another $\gamma_i$-perfect Hermitian form $\A'$ whose Voronoi domain contains $\mathcal{S}$. We certainly have $\V_i(\A)\cap \V_i(\A_\rho)=\mathcal{S}$. Because of $\V_i(\A') \supsetneq \mathcal{S}$ there is a common interior point of $\V_i(\A')$ and $\V_i(\A)$ or $\V_i(\A_\rho)$. We then obtain $\A'=\A$ or $\A'=\A_\rho$.\beweisende

\begin{definition}The $\gamma_i$-perfect Hermitian form $\A_\rho$ from the previous theorem is called the contiguous form to $\A$ (through the facet $\mathcal{S}$). \end{definition}

With regard to the implementation of the algorithm note that $\rho$ in Theorem \ref{contiguous} is a rational number, provided that $\min\nolimits_{L_i}(\A) \in \MQ$ and the associated facet vector is chosen from $K^{n\times n}$.

Now Theorem \ref{contiguous} constitutes a major part of the desired algorithm, provided one knows a $\gamma_i$-perfect form with which to use that theorem. The following lemma, which can be proven analogously to \ref{contiguous}, provides a way to find a first $\gamma_i$-perfect form.

\begin{lemma}\label{firstperfect}Let $\A\in\H_n$ be positive definite but not $\gamma_i$-perfect and $R\in \langle \V_i(\A)\rangle ^\perp$. Put $\A_t := \A + t\cdot R$ for $t\in \MR$. Then there is exactly one $\rho \in (0,\infty]$ such that $\min\nolimits _{L_i}(\A_t) = \min \nolimits _{L_i} (\A)$ for all $0\leq t \leq \rho$. Furthermore we have $\dim ( \langle \V_i(\A_\rho) \rangle ) > \dim ( \langle \V_i(\A) \rangle   )$.\\
If $t>\rho$ then $\A_t$ has smaller minimum than $\A$ or $\A_t$ is not positive definite anymore.\end{lemma}

This lemma shows that we may choose an arbitrary positive definite Hermitian form. If it is not $\gamma_i$-perfect we can apply the above process repeatedly until we obtain a perfect form, since the dimension of the subspace spanned by the Voronoi domain will strictly increase in every step.

\begin{definition} The Voronoi graph of the lattice $L_i$ is the graph whose vertices are the classes of $\gamma_i$-perfect Hermitian forms in $\MR_{>0}\backslash \H_n / \GL (L_i)$. Two vertices are connected by an edge whenever the respective classes contain contiguous perfect forms.\end{definition}
Counting the number of facets for each Voronoi domain, the Voronoi graph can be made into a weighted directed graph.

The following theorem shows that Theorem \ref{contiguous} can indeed be used to formulate an algorithm which produces a complete list of the $\gamma_i$-perfect forms.

\begin{thm}\label{graph}The Voronoi graph is a finite connected graph.\end{thm}

\textbf{Proof:} The finiteness is shown in Theorem \ref{finite}. Now let $\A$ and $\mathfrak{B}$ be two $\gamma_i$-perfect Hermitian forms of equal minimum. Choose an interior point $T\in \V_i(\mathfrak{B})$. If $T\in \V_i(\A)$, then $\A= \mathfrak{B}$ by Theorem \ref{interior}.\\
Otherwise there is a facet vector of $\A$ such that $\trace (TR) <0$. Define $\A_1$ to be the contiguous form to $\A$ through the facet defined by $R$. Then
\begin{equation*}\trace (\A_1 T) < \trace ( \A T) . \end{equation*}
Either now we have $T\in \V_i(\A_1)$ or we can repeat the above method, diminishing $\trace (\A_i T)$ for each new $\A_i$ obtained by this process. However, this process must be finite by \cite[7.3.2]{mar03}, so that there is indeed a finite sequence of contiguous forms $\A_0 = \A , ~ \A_1 , ... , \A_r = \mathfrak{B}$. \beweisende

\subsection{Implementation}

The algorithm can be implemented taking recourse to the algorithms implemented for $\MZ$-lattices in Magma \cite{magma} combined with the program QHull \cite{qhull}, which calculates the facets and facet vectors of the considered Voronoi domains. We may obtain the Gram matrix of a $\MZ$-lattice from a Hermitian form $\A$ with the so-called trace form.

\begin{definition}\label{def:traceform}Let $\A\in \H_n$ such that the entries of $\A$ are contained in $K$. $L_i$ is a free $\MZ$-module of rank $2n$ with basis $B:=(b_1,...,b_{2n})$. We then call
\begin{equation*}\frac{1}{2} ( \tr_{K/\MQ} ( b_i \A b_j^* ) ) _{i,j} \in \MQ^{2n \times 2n} \end{equation*}
the trace form of $\A$ with respect to the basis $B$. In this case $\tr_{K/\MQ}$ denotes the field trace of $K/\MQ$.\end{definition}

The following theorem is well-known.

\begin{thm}\label{isometry}Let $a, ~b ~:~ K^n \times K^n \to K$ be two non-degenerate Hermitian forms on $K^n$. For $\varphi ~:~ K^n \to K^n$ and all $v,u \in K^n$ the following are equivalent.
\begin{enumerate}
\item{$\varphi$ is $\MQ$-linear and $\tr_{K/\MQ}(a(\varphi(v),\varphi(u))) = \tr_{K/\MQ}(b(v,u))$ and $\tr_{K/\MQ}(\omega a(\varphi(v), \varphi(u))) = \tr_{K/\MQ}(\omega b(v,u))$}

\item{$\varphi$ is $K$-linear and $a(\varphi(v),\varphi(u))=b(v,u)$}
\end{enumerate}
\end{thm}

Thanks to this theorem we can test equivalence of Hermitian forms by testing isometry of $\MZ$-lattices using the trace form and \cite{ps97}.

We can now formulate the complete algorithm.

\begin{algorithm}[Voronoi] \label{algo}
\textbf{Input:} Imaginary quadratic number field $K$, dimension $n$, an ideal $\a_i$ in order to determine the lattice $\O_K^{n-1} \oplus \a_i$.\\
\textbf{Output:} A list $\mathcal{L}_i$ of all $\gamma_i$-perfect Hermitian forms.\\
\begin{enumerate}
\item{Find an integral representative $\a_i$ of minimal norm of the $i^{\text{th}}$ ideal class and a $\MZ$-basis of $L_i = \O_K^{n-1} \oplus \a_i$.}
\item{Use Lemma \ref{firstperfect} in order to determine a first $\gamma_i$-perfect form. Add it to $\mathcal{L}_i$ and to a list $\mathcal{T}_i$.}
\item{For every $T\in \mathcal{T}_i$ determine the contiguous forms using Theorem \ref{contiguous}.\\
Check the contiguous forms for equivalence with all forms in $\mathcal{L}_i$ using \ref{isometry}. Add new forms to $\mathcal{L}_i$ and $\mathcal{T}_i$. Contiguous forms which are equivalent to a Hermitian form in $\mathcal{L}_i$ may be ignored.\\
Put $\mathcal{T}_i := \mathcal{T}_i - \{ T \}$ and repeat this step until $\mathcal{T}_i = \emptyset$.}
\end{enumerate}
\end{algorithm}

\begin{thm}The above algorithm terminates and yields a complete list of all $\gamma_i$-perfect Hermitian forms. \end{thm}
\textbf{Proof:} This follows from Theorem \ref{graph}. The finiteness of the graph implies that the algorithm terminates. Since the graph is connected the list $\mathcal{L}_i$ of $\gamma_i$-perfect forms is complete. \beweisende

\begin{remark}Since the study of $n$-dimensional perfect Hermitian forms leads to the study of $2n$-dimensional perfect quadratic forms, the algorithm is, at present, only efficient up to dimension $3$.\\
In dimension $4$ we will encounter the $\MZ$-lattice $\mathds{E}_8$ \cite{hkn12}, whose Voronoi domain has 25075566937584 facets (see also \cite{dsv07}). Enumerating all perfect forms with the conventional Voronoi algorithm is therefore not feasible in dimensions greater or equal than $4$.\end{remark}

\section{An application of Voronoi's algorithm}\label{opgen}

In this section we present an example of how to apply the results from \cite{Opg01} in order to obtain a set of generators for the group $\GL(L)$, where $L$ is an $\O_K$-lattice. To our knowledge, this is the first example of a general procedure which allows for the computation of such a group, in particular when $L$ is non-free as a  $\O_K$-module (see below for an explicit example). More specifically, we make use of the following theorem, which is a reformulation of \cite[Theorem 2.2]{Opg01} in our context, and is essentially a consequence of Bass an Serre's theory \cite{Ser77}:

\begin{thm}\cite[Theorem 2.2]{Opg01}\label{opg} 
Let $\left\lbrace \a_1, \dots , \a_{h_K}\right\rbrace$ be a set of representatives of the class-group of $K$. For each $i=1, \dots, h_K$, we set $L_i=\O_K ^{n-1}\oplus \a_i$ and fix :
\begin{itemize}
\item a finite set $\mathcal{L}_i$ of representatives of the $\gamma_i$-perfect forms,
\item the finite set $\mathcal{M}_i$ of  $\gamma_i$-perfect forms not in $\mathcal{L}_i$ that are contiguous, in the Voronoï graph, to a form in $\mathcal{L}_i$,
\item for each $\mathcal{A} \in \mathcal{M}_i$, a matrix $U_{\mathcal{A}} \in \GL(L_i)$ such that $\mathcal{A}[U_A^{-1}] \in \mathcal{L}_i$.
\end{itemize} 
Then \begin{equation*}
\GL(L_i)=\left\langle \stab_{\GL(L_i)}(\mathcal P), ~U_{\mathcal{A}} ~|~\mathcal P \in \mathcal{L}_i, ~ \mathcal{A} \in \mathcal{M}_i \right\rangle 
\end{equation*}
\end{thm}

It should be noted that the finite groups $\stab_{\GL(L_i)}(\mathcal P)$ may of course be calculated without having to calculate the whole group $\GL(L_i)$. In order to do so, we use $\mathcal{P}$ to construct a $2n$-dimensional $\MZ$-lattice via the trace-form \ref{def:traceform}. We then use an algorithm which is based on \cite{ps97} and which is implemented in Magma\footnote{The Magma procedure may be called by \texttt{AutomorphismGroup}.} \cite{magma} to determine the finite matrix group which stabilizes the constructed lattice and the trace form matrix of $\omega \mathcal{P}$. From Theorem \ref{isometry} it follows that this yields a group of $\mathcal{O}_K$-linear automorphisms which stabilize $\mathcal{P}$, if we reduce the obtained $2n\times 2n$-matrices over $\MQ$ to $n\times n$-matrices over $K$ by comparing the action of the matrices on vector space bases of $K^n$ and $\MQ^{2n}$.

\begin{example}
Concretely, we deal with the general linear group of the non-free lattice $L_2:=\O_{\MQ(\sqrt{-15})}\oplus \left\langle 2, ~ \frac{\sqrt{-15}+1}{2}-1 \right\rangle$, see Section \ref{results} for a detailed exposition of the results of Voronoi's algorithm for this lattice.\\
Note that in the case of a free lattice $L$ there are other methods to obtain generators of $\GL(L)$, see \cite{EGM98, Bia92}. To our knowledge, no such computation has been done in the case of a non-free lattice. In addition, the approach presented here is easily adapted to lattices of higher dimension, the only restriction being the computational difficulty of studying the Voronoi domain and enumerating all neighbours of a given perfect form. Opgenorth's method may also be applied to all other lattices treated in Section \ref{results}.

In the present case there is just one perfect form, $\mathfrak{P}=\{\mathcal{P}\}$. For every $\gamma_2$-perfect form $\mathcal{A}\neq\mathcal{P}$ which is contiguous to $\mathcal{P}$, we determine $U_\A \in \GL(L_2)$ such that $\A=\mathcal{P}[U_\A]$. Then, according to theorem \ref{opg}, one has

$$\GL(L_2) = \langle \stab_{\GL(L_2)}(\mathcal{P}), ~ U_\A ~|~ \A \text{ contiguous to } \mathcal{P} \rangle.$$

In the concrete example we have {\tiny $\mathcal{P} = \begin{pmatrix} 1 & \frac{1}{10} ( 5+\sqrt{-15} ) \\ \frac{1}{10} (5-\sqrt{-15}) & \frac{1}{2} \end{pmatrix}$}. The stabiliser of $\mathcal{P}$ is isomorphic to $C_3 \rtimes C_4$ and it is generated by
{
\small
\begin{equation*}
\begin{pmatrix}1 & \frac{1}{2}(-3-\sqrt{-15}) \\  \frac{1}{4}(3-\sqrt{-15}) & -2\end{pmatrix} , \quad \begin{pmatrix}\frac{1}{2}(1-\sqrt{-15}) & \frac{1}{2}(-5+\sqrt{-15}) \\ \frac{1}{4}(-1-\sqrt{-15}) & \frac{1}{2}(-1+\sqrt{-15})\end{pmatrix}.
\end{equation*}
}

$\mathcal{P}$ has eight neighbours, the corresponding $U_\A$ can be chosen as
{
\small
\begin{equation*}\begin{split}
&\left(\begin{smallmatrix}-4-\sqrt{-15} & 2\sqrt{-15} \\ -4 & 4+\sqrt{-15}\end{smallmatrix} \right) , ~ \left( \begin{smallmatrix}4+\sqrt{-15} & \frac{1}{2}(-1-3\sqrt{-15}) \\ 3 & \frac{1}{2}(-5-\sqrt{-15})\end{smallmatrix}\right) , ~ \frac{1}{2} \left( \begin{smallmatrix}7-\sqrt{-15} & -15-\sqrt{-15} \\ 1-\sqrt{-15}  &   -7+\sqrt{-15} \end{smallmatrix}\right), ~ \left(\begin{smallmatrix}1 & \frac{1}{2}(-3-\sqrt{-15}) \\ 0 & -1 \end{smallmatrix} \right), \\
&\left(\begin{smallmatrix}2-\sqrt{-15} & \frac{1}{2}(-5+\sqrt{-15}) \\ \frac{1}{4}(-3-3\sqrt{-15}) & \frac{1}{2}(-1+\sqrt{-15})\end{smallmatrix}\right), ~ \frac{1}{2}\left( \begin{smallmatrix} -3+\sqrt{-15} & 5-\sqrt{-15} \\ -1+\sqrt{-15} & 3-\sqrt{-15} \end{smallmatrix} \right) , ~ \left(\begin{smallmatrix}1 & \frac{1}{2}(-7-\sqrt{-15}) \\ 0 & -1 \end{smallmatrix} \right) , ~ \left( \begin{smallmatrix}1&-2 \\ 0&-1 \end{smallmatrix}\right)
\end{split}
\end{equation*}
}

In conclusion, the ten matrices presented above constitute a set of generators of the group $\GL \left( L_2 \right)$. The matrices were obtained by explicit calculations with \cite{magma,qhull}.
\end{example}

\begin{remark} From this computation we also obtain an argument which shows that $\GL(L_2) \not \cong \GL_2(\O_K)$. Otherwise, $\GL_2(\O_K)$ would contain a subgroup isomorphic to $G:=\stab_{\GL(L_2)}(\mathcal{P})$, the lattere being isomorphic to $C_3 \rtimes C_4$. Consequently, $G$ would fix a free $\O_K$-lattice, but the $\O_K G$-lattices in $K^2$ are all of the form $\c L_2$ for some ideal $\c \trianglelefteq \O_K$. Note that the $\O_K$-order generated by the matrices in $G$ is of index $3$ in the maximal order $\mathrm{End}_{\O_K}(L_2)$. Nevertheless one computes (with Magma \cite{magma}) that the sublattices of $L_2$ of $3$-power index are the same as for $\mathrm{End}_{\O_K}(L_2)$. \end{remark}

\section{Computational results}\label{results}

\subsection{Dimension 2}

An implementation of Algorithm \ref{algo} in Magma \cite{magma} using QHull \cite{qhull} to treat the Voronoi domains produces the following computational results. The implementation is based on the implementation of Voronoi's algorithm in \cite{meyer} for the case $h_K =1$.

The following table presents some of the obtained results for the field $K= \MQ(\sqrt{-d})$ and the lattice $L=\O_K \oplus \a$.\\
The table also contains some important invariants, namely the determinant relative to $L$, the number of shortest vectors up to the relation $\sim$ as defined in Section 3, the number of facets of the Voronoi domain and the automorphism group of $(L,P)$, i.e. the group of $\O_K$-module automorphisms which fix $P$.\\
Note that all Hermitian forms are scaled such that their minimum is equal to $1$.

The occurring fields are ordered by the absolute values of their discriminants, since the calculated results suggest that the number of perfect Hermitian forms grows as the discriminant increases. However we do not know if this is always the case or why the number of perfect forms should increase with the discriminant.

{
\tiny
\begin{center}

\begin{tabular}{|c|c|c|c|c|c|c|}
\hline $d$ & $\a$ & $P$ & $\det_{L}(P)$ & $|S_{L}(P)|$ & facets & $\mathrm{Aut}(L,P)$ \\ 
\hline \hline $15$ & $\O_K$ & $\left( \begin{smallmatrix}1 & \frac{1}{6}(3+\sqrt{-15}) \\ \frac{1}{6}(3-\sqrt{-15}) & 1\end{smallmatrix}\right)$ & $\frac{1}{3}$ & $6$ & $8$ & $C_6$ \\
&& $\left( \begin{smallmatrix}2 & \frac{1}{10}(15+3\sqrt{-15}) \\ \frac{1}{10}(15-3\sqrt{-15}) & 2\end{smallmatrix} \right)$ & $\frac{2}{5}$ & $4$ & $4$ & $C_4$ \\
\hline & $\left \langle 2, ~ \frac{\sqrt{-15}+1}{2}-1 \right \rangle $ & $\left(\begin{smallmatrix}1 & \frac{1}{10}(5+\sqrt{-15}) \\ \frac{1}{10}(5-\sqrt{-15}) & \frac{1}{2}\end{smallmatrix}\right)$ & $\frac{1}{5}$ & $12$ & $8$ & $C_3 \rtimes C_4$ \\
\hline \hline $5$ & $\O_K$ & $\left( \begin{smallmatrix} 1 & \frac{1}{10} (5+3\sqrt{-5}) \\ \frac{1}{10}(5-3\sqrt{-5}) & 1 \end{smallmatrix} \right)$ & $\frac{3}{10}$ & $6$ & $5$ & $C_6$ \\
 && $\left( \begin{smallmatrix}1 & \frac{1}{5} ( 5+2\sqrt{-5}) \\ \frac{1}{5}(5-2\sqrt{-5}) & 2 \end{smallmatrix} \right)$ & $\frac{1}{5}$ & $8$ & $6$ & $Q_8$ \\
\hline & $\langle 2, ~ 1+\sqrt{-5} \rangle$ & $\left( \begin{smallmatrix} 1 & \frac{1}{10}(5+2\sqrt{-5}) \\ \frac{1}{10} (5-2\sqrt{-5}) & \frac{1}{2} \end{smallmatrix} \right)$ & $\frac{1}{10}$ & $24$ & $14$ & $\mathrm{SL}(2,3)$ \\
\hline \hline $23$ & $\O_K$ & $\left(\begin{smallmatrix}1 & \frac{1}{46}(23+7\sqrt{-23}) \\ \frac{1}{46}(23-7\sqrt{-23}) & 1\end{smallmatrix}\right)$ & $\frac{5}{23}$ & $9$ & $8$ & $C_6$ \\
&& $\left(\begin{smallmatrix}3 & \frac{1}{46}(115+15\sqrt{-23}) \\ \frac{1}{46}(115-15\sqrt{-23}) & 3\end{smallmatrix} \right)$ & $\frac{7}{23}$ & $6$ & $5$ & $C_4$ \\
\hline \hline $6$ & $\O_K$ & $\left( \begin{smallmatrix} 1 & \frac{1}{6} (3+2\sqrt{-6}) \\ \frac{1}{6}(3-2\sqrt{-6}) & 1 \end{smallmatrix} \right)$ & $\frac{1}{12}$ & $24$ & $26$ & $\mathrm{SL}(2,3)$ \\
\hline & $\langle 2, ~ \sqrt{-6} \rangle$ & $\left( \begin{smallmatrix} 1 & \frac{1}{6}(3+\sqrt{-6}) \\ \frac{1}{6} (3-\sqrt{-6}) & \frac{1}{2} \end{smallmatrix} \right)$ & $\frac{1}{6}$ & $8$ & $6$ & $Q_8$ \\
 & & $\left( \begin{smallmatrix} 1 & \frac{1}{4}(3+\sqrt{-6}) \\ \frac{1}{4} (3-\sqrt{-6}) & 1 \end{smallmatrix}\right)$ & $\frac{1}{8}$ & $12$ & $8$ & $C_3 \rtimes C_4$ \\
\hline
\end{tabular}\\
Table 1: Computational results (1/2)
\end{center}
}

{
\tiny
\begin{center}
\begin{tabular}{|c|c|c|c|c|c|c|}
\hline $10$ & $\O_K$ & $\left( \begin{smallmatrix} 1 & \frac{1}{4} (2+\sqrt{-10}) \\ \frac{1}{4}(2-\sqrt{-10}) & 1 \end{smallmatrix} \right)$ & $\frac{1}{8}$ & $6$ & $8$ & $C_6$ \\
&& $\left( \begin{smallmatrix}1 & \frac{1}{30}(15+11\sqrt{-10}) \\ \frac{1}{30}(15-11\sqrt{-10}) & \frac{5}{3} \end{smallmatrix} \right)$ & $\frac{13}{180}$ & $12$ & $10$ & $C_4$ \\
\hline & $\langle 2, \sqrt{-10} \rangle$ & $\left(\begin{smallmatrix} 1 & \frac{1}{20}(10+3\sqrt{-10}) \\ \frac{1}{20} (10-3\sqrt{-10}) & \frac{1}{2} \end{smallmatrix}\right)$ & $\frac{1}{20}$ & $24$ & $14$ & $\mathrm{SL}(2,3)$ \\
&& $\left(\begin{smallmatrix} 4 & \frac{1}{20}(45+14\sqrt{-10}) \\ \frac{1}{20} (45-14\sqrt{-10}) & \frac{5}{2} \end{smallmatrix}\right)$ & $\frac{3}{40}$ & $12$ & $8$ & $C_3 \rtimes C_4$ \\
&& $\left(\begin{smallmatrix}5 & \frac{1}{20}(55+14\sqrt{-10}) \\ \frac{1}{20}(55-14\sqrt{-10}) & \frac{5}{2}\end{smallmatrix}\right)$ & $\frac{3}{40}$ & $12$ & $8$ & $C_3\rtimes C_4$ \\
\hline \hline $21$ & $\O_K$ & $\left( \begin{smallmatrix} 1 & \frac{1}{6} (3+\sqrt{-21}) \\ \frac{1}{6} (3-\sqrt{-21}) & 1 \end{smallmatrix} \right)$ & $\frac{1}{6}$ & $6$ & $8$ & $C_6$ \\
&& $\left( \begin{smallmatrix}1 & \frac{1}{14} (7+4\sqrt{-21}) \\ \frac{1}{14} ( 7 - 4\sqrt{-21}) & 2 \end{smallmatrix} \right)$ & $\frac{1}{28}$ & $24$ & $26$ & $C_3 \rtimes C_4$ \\
&& $\left(\begin{smallmatrix}1 & \frac{1}{14}(14+3\sqrt{-21}) \\ \frac{1}{14}(14-3\sqrt{-21}) & 2 \end{smallmatrix} \right)$ & $\frac{1}{28}$ & $24$ & $26$ & $C_3 \rtimes C_4$ \\
&& $\left( \begin{smallmatrix}13 & \frac{1}{6}(39+23\sqrt{-21}) \\ \frac{1}{6}(39-23\sqrt{-21}) & 27 \end{smallmatrix} \right)$ & $\frac{1}{6}$ & $6$ & $8$ & $C_6$ \\
&& $\left( \begin{smallmatrix}5 & \frac{1}{14}(35+19\sqrt{-21}) \\ \frac{1}{14}(35-19\sqrt{-21}) & 9\end{smallmatrix} \right)$ & $\frac{1}{14}$ & $12$ & $8$ & $C_6$ \\
&& $\left( \begin{smallmatrix}23 & \frac{1}{7}(84+45\sqrt{-21}) \\ \frac{1}{7}(84-45\sqrt{-21}) & 44 \end{smallmatrix} \right)$ & $\frac{1}{7}$ & $8$ & $6$ & $C_4$ \\
\hline & $\langle 2 , ~ \sqrt{-21}-1 \rangle$ & $\left( \begin{smallmatrix}1 & \frac{1}{42}(21+4\sqrt{-21}) \\ \frac{1}{42}(21-4\sqrt{-21}) & \frac{1}{2} \end{smallmatrix} \right)$ & $\frac{5}{42}$ & $8$ & $6$ & $C_4$ \\
&& $\left(\begin{smallmatrix}2 & \frac{1}{14}(14+3\sqrt{-21}) \\ \frac{1}{14}(14-3\sqrt{-21}) & 1\end{smallmatrix}\right)$ & $\frac{1}{14}$ & $8$ & $6$ & $C_4$ \\
&& $\left(\begin{smallmatrix}\frac{7}{3} & \frac{1}{42} (35+6\sqrt{-21}) \\ \frac{1}{42} ( 35-6\sqrt{-21} ) & \frac{1}{2}\end{smallmatrix} \right)$ & $\frac{11}{126}$ & $8$ & $6$ & $C_2$ \\
&& $\left(\begin{smallmatrix}\frac{7}{3} & \frac{1}{14} (21+2\sqrt{-21}) \\ \frac{1}{14}(21-2\sqrt{-21}) & \frac{7}{6} \end{smallmatrix} \right)$ & $\frac{11}{126}$ & $8$ & $6$ & $C_2$ \\
&& $\left(\begin{smallmatrix}7 & \frac{1}{42} ( 147 + 32\sqrt{-21}) \\ \frac{1}{42}(147-32\sqrt{-21}) & \frac{7}{2} \end{smallmatrix}\right)$ & $\frac{5}{42}$ & $8$ & $6$ & $C_4$ \\
\hline & $\langle 5, ~ \sqrt{-21}-2 \rangle$ & $\left( \begin{smallmatrix}1 & \frac{1}{105}(42+4\sqrt{-21}) \\ \frac{1}{105}(42-4\sqrt{-21}) & \frac{1}{5} \end{smallmatrix} \right)$ & $\frac{1}{21}$ & $16$ & $10$ & $Q_8$ \\
&& $\left(\begin{smallmatrix}25 & \frac{1}{105}(1092+94\sqrt{-21}) \\ \frac{1}{105}(1092-94\sqrt{-21}) & 5\end{smallmatrix} \right)$ & $\frac{1}{21}$ & $16$ & $10$ & $Q_8$ \\
&& $\left(\begin{smallmatrix} 10 & \frac{1}{30}(129+8\sqrt{-21}) \\ \frac{1}{30}(129-8\sqrt{-21}) & 2\end{smallmatrix} \right)$ & $\frac{1}{12}$ & $8$ & $6$ & $C_4$ \\
&& $\left(\begin{smallmatrix}13 & \frac{1}{30}(156+17\sqrt{-21}) \\ \frac{1}{30}(156-17\sqrt{-21}) & \frac{13}{5}\end{smallmatrix}\right)$ & $\frac{1}{12}$ & $8$ & $6$ & $C_4$ \\
&& $\left(\begin{smallmatrix}\frac{4}{3} & \frac{1}{45}(21+2\sqrt{-21}) \\ \frac{1}{45}(21-2\sqrt{-21}) & \frac{1}{5}\end{smallmatrix} \right)$ & $\frac{1}{27}$ & $16$ & $10$ & $C_4$ \\
\hline & $\langle 3,  ~\sqrt{-21} \rangle$ & $\left(\begin{smallmatrix}1 & \frac{1}{18}(9+\sqrt{-21}) \\ \frac{1}{18}(9-\sqrt{-21}) & \frac{1}{3}\end{smallmatrix} \right)$ & $\frac{1}{18}$ & $12$ & $8$ & $C_6$ \\
&& $\left(\begin{smallmatrix}2 & \frac{1}{12}(11+\sqrt{-21}) \\ \frac{1}{12}(11-\sqrt{-21}) & \frac{1}{2}\end{smallmatrix}\right)$ & $\frac{1}{24}$ & $16$ & $10$ & $C_4$\\ 
&& $\left(\begin{smallmatrix}11 & \frac{1}{42}(217+23\sqrt{-21}) \\ \frac{1}{42}(217-23\sqrt{-21}) & 3\end{smallmatrix} \right)$ & $\frac{1}{42}$ & $48$ & $26$ & $\mathrm{SL}(2,3)$ \\
&& $\left(\begin{smallmatrix}5 & \sqrt{1}{42}(91+5\sqrt{-21}) \\ \frac{1}{42}(91-5\sqrt{-21}) & 1\end{smallmatrix}\right)$ & $\frac{1}{42}$ & $48$ & $26$ & $\mathrm{SL}(2,3)$ \\

\hline
\end{tabular}\\
Table 2: Computational results (2/2)
\end{center}
}

Note that in the case of $K=\MQ(\sqrt{-23})$ we have $\cl_K \cong C_3$ and therefore $\cl_K / \cl_K^2 \cong \{1\}$. Therefore, by virtue of \ref{simplification}, we merely consider the free lattice $\O_K^2$.

For the above mentioned fields we obtain the following Hermite constants.

\begin{center}
\begin{tabular}{|c||c|c|c|c|c|c|}
\hline $d$ & $15$ & $5$ & $23$ & $6$ & $10$ & $21$ \\
\hline $\gamma_{2,\MQ(\sqrt{-d})}$ & $\sqrt{5}$ & $\sqrt{10}$ & $\sqrt{\frac{23}{5}}$ & $\sqrt{12}$ & $\sqrt{20}$ & $\sqrt{42}$ \\
\hline
\end{tabular}	
\end{center}

In the case of $\gamma_{2,\MQ(\sqrt{-21})}$ we observe the existence of two inequivalent perfect Hermitian forms over the lattice $\O_{\MQ(\sqrt{-21})} \oplus \langle 3, ~ \sqrt{-21} \rangle$ which both realise the maximum value $\sqrt{42}$.

For illustration, we present the Voronoi graphs for the field $K=\MQ(\sqrt{-10})$. In the left one the underlying lattice is $\O_K^2$, in the right one we work over $\O_K \oplus \langle 2 , ~ \sqrt{-10} \rangle$.\\
The perfect forms in the graph are numbered in the order in which they appear in the above table. The marked vertex denotes the perfect form which realises the global maximum of Hermite's function. The arrows originating at a vertex $P$ point to the perfect neighbours of $P$, the weight is the number of facets of $\V(P)$ through which $P$ and the respective perfect neighbour are contiguous.

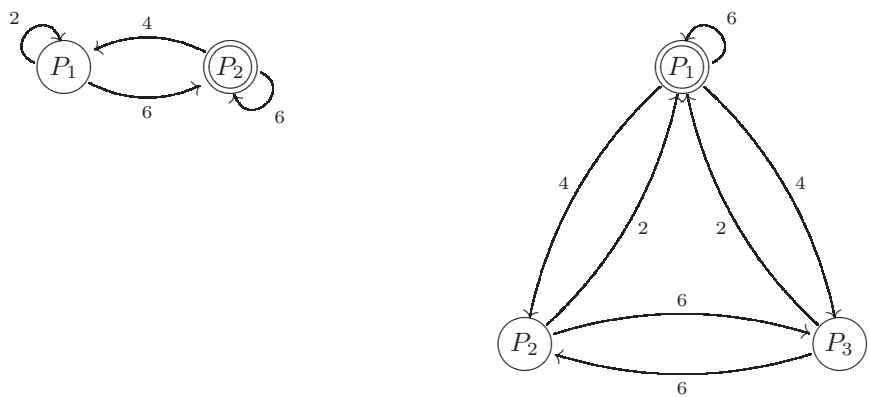
\begin{center}
{
\small
\begin{tabular}{ccc}
$$
\entrymodifiers={++[o][F-]}
\SelectTips{cm}{}
\xymatrix @-1pc { P_{1} \ar@(l,u)^2 \ar@/_ .4 cm/[rr]_6 & *++{~} & *++[o][F=]{P_{2}} \ar@(r,d)^6 \ar@/_ .4 cm/[ll]_4
}

&\phantom{\hspace{1.8cm}} &
\entrymodifiers={++[o][F-]}
\SelectTips{cm}{}
\xymatrix @-1pc { *++{}& *++{} & *++[o][F=]{P_{1}} \ar@(r,u)_6 \ar@/_.4cm/[ddddll]_4 \ar@/^.4cm/[ddddrr]^4 & *++{} & *++{} \\
*++{} & *++{} & *++{~} & *++{} & *++{} \\
*++{} & *++{} & *++{~} & *++{} & *++{} \\
*++{} & *++{} & *++{~} & *++{} & *++{} \\
P_{2} \ar@/_.4cm/[uuuurr]_2 \ar@/^.4cm/[rrrr]^6 & *++{} & *++{} & *++{} & P_{3} \ar@/^.4cm/[llll]^6 \ar@/^.4cm/[uuuull]^2
}
$$
\end{tabular}
}
\end{center}

\subsection{Dimension 3}

This section is devoted to the results obtained in dimension 3. In this situation there are too many perfect Hermitian forms to list all Gram matrices in print. We therefore restrict ourselves to giving the numbers of perfect forms. Gram matrices can be obtained from \url{http://www.math.rwth-aachen.de/~Oliver.Braun/perfect.html}. We again rely on \ref{simplification} in order to omit some cases.

Note that the considered fields are of class number 2. Since $\cl_K / \cl_K^3 \cong \{1\}$, we only have to consider the free lattice $\O_K^3$, which led to the following results for the fields $\MQ(\sqrt{-d})$.

{
\begin{center}
\begin{tabular}{|c||c|c|c||c|c|}
\hline $d$ & $15$ & $5$ & $6$  & $10$ & $13$ \\
\hline number of perfect forms & $11$ & $92$ & $271$ & $4236$ & $\geq 2746$  \\
\hline $\gamma_{3,\MQ(\sqrt{-d})}$ & $\sqrt[3]{15}$ & $\sqrt[3]{20}$ & $\sqrt[3]{24}$ &  $\sqrt[3]{\frac{60835}{911}}$ & $\geq \sqrt[3]{ \frac{632684}{6749}}$ \\
\hline
\end{tabular}
\end{center}
}

Over the field $\MQ(\sqrt{-5})$ there are two inequivalent Hermitian forms which realise the global maximum of Hermite's function, over $\MQ(\sqrt{-6})$ there are seven. In the case of $\MQ(\sqrt{-15})$ the maximum is unique.

 \providecommand{\bysame}{\leavevmode\hbox to3em{\hrulefill}\thinspace}
\providecommand{\MR}{\relax\ifhmode\unskip\space\fi MR }
\providecommand{\MRhref}[2]{%
  \href{http://www.ams.org/mathscinet-getitem?mr=#1}{#2}
}
\providecommand{\href}[2]{#2}

\end{document}